\documentclass[centertags,leqno,11pt]{article}

\usepackage{amsmath}
\usepackage{amssymb}
\usepackage{latexsym}
\usepackage{amsxtra}
\usepackage{amscd}
\usepackage{theorem}

\usepackage{graphics}
\usepackage{epic}

\theoremstyle{change}

{\theorembodyfont{\slshape}

\newtheorem{thm}{Theorem.}[section]

\newtheorem{lem}[thm]{Lemma.}
\newtheorem{prop}[thm]{Proposition.}

\newtheorem{defn}[thm]{Definition.}}

{\theorembodyfont{\rmfamily}

\newtheorem{rem}[thm]{Remark.}

}

\renewcommand{\em}{\sl}

\newcommand{\proof}{\noindent {\bf Proof:\ }}
\newcommand{\Endproof}{\hspace*{\fill} $\Box$ \vspace{1ex} \noindent }

\makeatletter
\renewcommand{\subsection}{\@startsection{subsection}{2}%
{\z@}{-3.25ex plus -1ex minus-.2ex}{-1em}{\bf}} \makeatother
\newcommand{\wV}{{\sqrt{5}}}
\newcommand{\wIII}{{\sqrt{3}}}

\newcommand{\PP}{\mathbb{P}}
\newcommand{\ZZ}{\mathbb{Z}}
\newcommand{\CC}{\mathbb{C}}
\newcommand{\RR}{\mathbb{R}}
\newcommand{\QQ}{\mathbb{Q}}
\newcommand{\NN}{\mathbb{N}}

\renewcommand{\AA}{\mathbb{A}}
\newcommand{\GG}{\mathbb{G}}

\newcommand{\V}{\mathcal{V}}

\newcommand{\K}{\mathcal{K}}

\newcommand{\p}{\mathfrak{p}}
\newcommand{\GL}{{\rm GL}}
\newcommand{\SL}{{\rm SL}}
\newcommand{\PGL}{{\rm PGL}}

\newcommand{\Gal}{{\rm Gal}}

\renewcommand{\c}{{\rm c}}

\newcommand{\MC}{{\rm MC}}

\newcommand{\A}{{A}}

\newcommand{\rk}{{\rm rk}}
\newcommand{\id}{{\rm id}}

\newcommand{\GGG}{{\cal G}}

\renewcommand{\L}{{\cal L}}

\newcommand{\M}{{ M}}

\newcommand{\mc}{{\rm mc}}

\numberwithin{equation}{subsection}
\numberwithin{thm}{subsection}
\theoremstyle{plain}

\begin{document}

\title{On globally nilpotent differential equations}

\author{
   Michael Dettweiler and 
 Stefan Reiter
  }

\maketitle 
\begin{abstract} In a previous work of the authors, a middle 
convolution operation on the category of 
 Fuchsian differential systems was introduced. In this note we show 
that the middle convolution of Fuchsian systems preserves the property of
global nilpotence. 

This leads to a globally nilpotent Fuchsian system 
of rank two which does not belong
to the known classes of globally nilpotent 
rank two systems. 
Moreover, we give 
a globally nilpotent Fuchsian system of 
rank seven whose differential Galois group is isomorphic
to the exceptional simple algebraic group of type $G_2.$

\end{abstract}


\section*{Introduction}\label{Introduction}

A unifying description of all  irreducible and physically rigid local
systems on the punctured affine line was given by Katz
\cite{Katz96}. The main tool therefore is a middle convolution
functor on the category of  perverse sheaves (loc.~cit., Chap.~5).
In \cite{DR}, the authors give a  purely algebraic analogon
of this convolution functor. 
This  functor is a functor of the category
 of finite dimensional 
$K$-modules of the free group $F_r$ on $r$
generators to itself ($K$ denoting a field). It depends on a scalar
$\lambda\in K^\times$ and is denoted by $\MC_\lambda.$ 

By the Riemann Hilbert correspondence (see \cite{Deligne70}), a construction parallel
to  $\MC_\lambda$ should exist in the category of Fuchsian systems
of differential equations.
In \cite{DR03}, such a  construction is given, leading to 
a  description of rigid Fuchsian systems which is parallel 
to Katz description of rigid local systems. The convolution 
depends on a parameter $\mu \in \CC$ and 
carries a Fuchsian system $F$ to another Fuchsian system, denoted 
by $\mc_\mu(F),$  see Section~\ref{sec13}.

In this note we study how the $\p$-curvature (for $\p$ a prime of a number field $K$)
  of a Fuchsian 
system $F$  having coefficients in the function field
$K(t)$
changes under the convolution process.  The $\p$-curvature is 
a matrix $\bar{C}_\p(F)$  with coefficients in the function 
field over a finite field which is obtained 
from a $p$-fold iteration of $F$ (where $p$ is the  prime number 
below $\p$) and reduction modulo $\p,$
see Section \ref{pcurv}.  
The $\p$-curvature matrices  encode many  arithmetic and geometric 
properties of
a Fuchsian system. For example, the Bombieri-Dwork conjecture 
predicts that  if the $\p$-curvature 
$\bar{C}_\p(F)$ is nilpotent for almost all 
primes $\p$ of $K$ (i.e., $F$ is {\em globally nilpotent}),
then $F$ is arising from geometry, see \cite{Andre89}. 
In this note we prove the following result (see Thm. \ref{cor21}):\\

\noindent {\bf Theorem 1:} {\em  Let 
$F$ be a Fuchsian system, let $\mu \in \QQ$ and let 
$\mc_\mu(F)$ be the middle convolution of $F$ with respect to $\mu.$   
Then the following holds:
\begin{enumerate}
\item
 If  
 ${\rm ord}_p(\mu) \geq 0$ and if the $p$-curvature of $F$ 
 is nilpotent of rank $k,$ then 
the $p$-curvature of the middle convolution 
$\mc_\mu(F)$ is nilpotent of rank $r\in \{k-1,k,k+1\}.$
\item If $F$ is globally nilpotent, then $\mc_\mu(F)$ is globally 
nilpotent.
\end{enumerate}}

Of course, the second statement of the theorem 
follows immediately from the first. 
The second statement can be deduced alternatively from the stability
of global nilpotence under pullback, tensor product, higher direct image,
see Katz \cite{Katznilpo}, Section 5.7 through 5.10 
(it follows from \cite{DR03}, Thm. 1.2,
and \cite{DeHabil}, Rem. 3.3.7, that $\mc_\mu$ 
corresponds under the Riemann-Hilbert
correspondence to the middle convolution $\MC_\chi$ of local systems $\V$
with Kummer sheaves - which is, by construction, a higher direct image sheaf).

The proof of the first statement 
 relies on  the closed formula of the 
$p$-curvature of Okubo systems (compare to Remark \ref{remput})
and on the fact that the middle convolution of a Fuchsian system 
is a factor system of an Okubo system.  \\

In the last section, we apply the above result in order to obtain two new 
 globally nilpotent differential equations:

In Section \ref{sec41}, we apply the middle convolution to 
the globally nilpotent 
Fuchsian system of rank two which appears in the work of Krammer
 \cite{Krammer}. This leads again to a globally 
nilpotent Fuchsian 
system of rank two.
It is a new type of a globally nilpotent 
rank two system because it is neither 
a pullback of a hypergeometric system, nor a
system associated to periods of Shimura curves
(see Thm. \ref{thm41} and the discussion in \cite{Krammer}, Section 11). 

In Section \ref{sec42} we consider  a 
globally nilpotent Fuchsian system $G$ of rank seven which is constructed by 
an sixfold iteration of middle convolutions and tensor products. 
Using the Riemann-Hilbert correspondence for $G$ (see 
Thm.~\ref{monodromy}), one can see that the Zariski 
closure of 
monodromy group  of $G$ is equal to the simple exceptional algebraic
group of type $G_2.$  The Fuchsian 
system $G$ can be viewed as the de Rham version of 
recent results of the authors on rigid local systems and motives with 
Galois group $G_2$ (\cite{DR06}). The motivic 
interpretation of the rigid local system 
which is given there, and  the Riemann-Hilbert 
correspondence (see Thm. \ref{monodromy}),
imply that  the Bombieri-Dwork conjecture
holds for the Fuchsian system~$G.$

The authors thank J. Aidan for pointing out a misprint in the differential
 equation which appears in Thm.~\ref{thm41} (in an earlier version of this paper) and
 N. Katz for valuable comments on the nilpotence of
 higher direct image connections.

\section{The Riemann-Hilbert correspondence of the middle convolution }\label{sec2}

\subsection{The 
tuple transformation $\MC_\lambda.$}
 Let $K$ be
 a field, let $V$ be a finite dimensional vectorspace over $K$ and let 
 $\M=(M_1,\ldots,M_r)$ be  an element of
$\GL(V)^r.$ 
For any  $\lambda \in K^\times$ one can construct another tuple 
of matrices 
$(N_1,\ldots,N_r)\in \GL(V^r)^r,$ as follows: For
$k=1,\ldots,r,$ $N_k$ maps a vector $(v_1,\ldots,v_r)^{\rm tr}$
$\in V^r$ to
\[ \left( \begin{array}{ccccccccc}
                  1 & 0 &  & \ldots& & 0\\
                   & \ddots &  & & &\\
                    & & 1 &&&\\
               (M_1-1) & \ldots&  (M_{k-1}-1)  & \lambda M_{k} &\lambda (M_{k+1}-1) & \ldots
&  \lambda (M_r-1) \\
     &&&&1&&\\
               &   &  & && \ddots  &   \\
                   0 &  &  & \ldots& &0 & 1
          \end{array} \right)\left(\begin{array}{c}
v_1\\
\vdots\\
\\\vdots\\
\\\vdots\\
\\ v_r\end{array}\right)
.\]

There are the following $\langle N_1,\ldots,N_r \rangle$-invariant
subspaces of  $V^r:$

\[ \K_k = \left( \begin{array}{c}
          0 \\
          \vdots \\
          0 \\
          \ker(M_k-1) \\
            0\\
           \vdots \\
           0
        \end{array} \right)  \quad \mbox{({\it k}-th entry)},\, k=1,\dots,r,\]
and
\[     \L=\cap_{k=1}^r \ker (N_k-1)={\rm ker}(N_r\cdots N_1 - 1).
\]
Let $\K:=\oplus_{i=1}^r\K_i.$

\begin{defn}{\rm Let
$\MC_\lambda(\M):=(\tilde{N}_1,\dots,\tilde{N}_r)\in
\GL(V^r/(\K+\L))^r,$ where $\tilde{N}_k$ is induced by the action
of $N_k$ on $V^r/(\K+\L).$ We call $\MC_\lambda(\M)$ the 
{\em  middle convolution} of $\M$ with
$\lambda.$}
\end{defn}

\subsection{The middle convolution of Fuchsian systems}\label{sec13}
Let
 $\A=(A_1,\ldots,A_r),\, A_k \in \CC^{n \times n}.$
 For $\mu \in \CC$ one can define
blockmatrices $B_k,\, k=1,\ldots,r,$ as follows:
  \[ B_k:= \left( \begin{array}{ccccccc}
                   0 &  &  & \ldots& & 0\\
                   & \ddots &  & & &\\
               A_1 & \ldots\,\, A_{k-1}&  A_k+ \mu  &  A_{k+1} & \ldots
&  A_r \\
               &   &  &  \ddots & &   \\
                   0 &    & \ldots& &  & 0
          \end{array} \right) \in \CC^{nr \times nr},\]
where $B_k$ is  zero outside the $k$-th block row.\\

The 
tuple 
\begin{equation}\label{eq133} c_\mu(\A):=(B_1,\ldots,B_r)\end{equation}
is called the {\em naive convolution} of $\A$ with $\mu.$ 
There are  the following left-$\langle B_1,\ldots,B_r
\rangle$-invariant subspaces of the column vector space $\CC^{nr}$
(with the tautological action of
$\langle B_1,\ldots,B_r \rangle$):\\

\[ {\frak{k}}_k = \left( \begin{array}{c}
          0 \\
          \vdots \\
            0\\
          \ker(A_k) \\
           0 \\
           \vdots \\
           0
        \end{array} \right) \quad \mbox{($k$-th entry)},\, k=1,\dots,r,\]
and
\[     {\frak{l}}=\cap_{k=1}^r {\rm ker}(B_k)
={\rm ker}(B_1+\ldots+B_r).
\]
Let ${\frak{k}}:=\oplus_{k=1}^r {\frak{k}}_k$ and 
 fix an isomorphism between $\CC^{nr}/({\frak{k}}+{\frak{l}})$
and $\CC^m.$

\begin{defn} {\rm   The tuple of matrices
$\mc_{\mu}(\A):=(\tilde{B}_1,\dots,\tilde{B}_r)\in \CC^{m\times m}, $
where $\tilde{B}_i$ is induced by the action of $B_i$ on
$\CC^m(\simeq \CC^{nr}/({\frak{k}}+{\frak{l}})),$ is called the {\em
 middle convolution} of $\A$ with $\mu.$}
\end{defn}

Let $S:=\CC\setminus \{t_1,\ldots,t_r\}$
 and $A:=(A_1,\ldots,A_r),\, A_i \in \CC^{n\times n}.$ The Fuchsian
system
$$F:Y'=\sum_{i=1}^r \frac{A_i}{t-t_i} Y$$ is called
 {\em the Fuchsian system associated to the tuple 
${A}.$}

\begin{defn}{\rm 
 Let $\A:=(A_1,\ldots,A_r),\, A_i \in \CC^{n\times n},$
and $\mu \in \CC.$ Let $F$ be the Fuchsian system associated to $A.$
Then the  Fuchsian system $F$
which is associated to the middle convolution tuple 
$\mc_\mu(\A)$
is denoted by $\mc_\mu(F)$ and is called  the {\em  middle
convolution} of $F$ with the parameter $\mu.$ The Fuchsian system 
which is associated to the naive convolution tuple $\c_\mu(\A)$
is denoted by $\c_\mu(F)$ and is called  the {\em  naive
convolution} of $F$ with the parameter $\mu.$}
\end{defn}

Fix a set of homotopy generators $\delta_i,\,i=1,\ldots,r,$
           of $\pi_1(S,t_0)$ by traveling from the base point
$t_0$ to $t_i,$ then encircling $t_i$ counterclockwise, and then 
going back to $t_0.$ Then, the analytic continuation of solutions
along $\delta_i$ defines a linear isomorphism $M_i,\,i=1,\ldots,r,$ 
of the vectorspace $V\simeq \CC^n$ of 
holomorphic solutions of $F$ which are defined  
in a small neighborhood of $t_0.$ 
We call the tuple 
$ (M_1,\ldots,M_r)$ the {\it monodromy tuple} of $F.$ The Riemann-Hilbert
correspondence says that $F$ is determined up to isomorphism
by $M,$ see \cite{Deligne70}. \\

The following result is an explicit realization
of the Riemann-Hilbert correspondence for a
convoluted Fuchsian system,  see \cite{DR03}, Thm. 6.8:

\begin{thm}\label{monodromy} 
Let $F$ be an irreducible Fuchsian system associated to 
$\A=(A_1,\ldots,A_r)\in (\CC^{n\times n})^r.$
 Assume that there exist two different
elements such that the local monodromy of $\V$ is nontrivial.  
Fix a set of homotopy generators 
$$\delta_1,\ldots,\delta_r\in 
\pi_1(\AA^1\setminus
\{t_1,\ldots,t_r\})$$ of the fundamental group 
as above. 
Let  
$M=(M_1,\ldots,M_r)\in \GL_n(\CC)^r$
be the monodromy tuple of  $F$ (with respect to $\delta_1,\ldots,\delta_r$).
 Assume that
$$\rk(A_i) = \rk(M_i-1),\,\,\,\,
\rk(A_1+\cdots + A_r+\mu)=\rk( \lambda\cdot M_1\cdots M_r-1).$$
Then the  monodromy tuple
of  $\mc_\mu(F)$ is given by $\MC_\lambda(M),$ where 
 $\lambda=e^{2\pi i \mu}.$ 
\end{thm}

\section{Transformation of the $p$-curvature under
$\mc_\mu$}\label{pcurv}

Let $K$ be a number field and let 
$F:Y'=CY$ be a system of linear differential equations, 
where $C\in
K(t)^{n\times n}.$ Successive application of differentiation
yields differential systems
\begin{equation} 
\label{eq0001}
F^{(n)}: Y^{(n)}={C_n}Y\quad {\rm for\,\, each }\quad n\in \NN_{>0}.\end{equation}
In the following, $\p$ always denotes a prime of $K$ which lies
over a prime number $p.$ For almost all primes $p,$   one can reduce
the entries of the matrices ${C_p}$ modulo $\p$
 in order
 to obtain the {\em $\p$-curvature matrices of $F$}
$$\bar{C}_\p=\bar{C}_\p(F):=C_p \mod \p.$$

\begin{defn}{\rm 
A  system of linear differential equations  $F$
{\em can be written in Okubo normal form}, if $F$ can be written as 
$$ F:Y'=(t-T)^{-1}BY\,,$$ where $B\in  \CC^{n\times
n}$ and $T$ is a diagonal matrix
 $T={\rm diag}(t_1,\ldots,t_n)$ with $ t_i \in \CC$
(here possibly $t_i=t_j$ for $i\not= j$). }
\end{defn}

The following proposition is obvious from the definitions:

\begin{prop}\label{proptriv} If  
$F$ is a Fuchsian system,  then
 the naive convolution $\c_\mu(F)$ of $F$ 
can be written   in Okubo normal form.
\end{prop}

An induction  yields the following formula for the 
$\p$-curvature matrix of a system in Okubo normal form:

\begin{lem}\label{recu} Let $$ F:Y'=CY=(t-T)^{-1}BY$$
be   a system of linear differential equations which can be written 
in Okubo normal form. 
 Then
 \begin{eqnarray}\label{eq2111}
{C_n}&=&(t-T)^{-1}(B-n+1)\cdot (t-T)^{-1}(B-n+2) \,\cdots \\
&& \quad \cdots \, (t-T)^{-1}
(B-1)\cdot (t-T)^{-1}B. \nonumber \end{eqnarray} Especially, if the matrix $B$ has coefficients in 
a number field $K$, then the $\p$-curvature ${\bar{C}_\p}$
of $F$ has the form
 $$(t-T)^{-1}(B-p+1)\cdot (t-T)^{-1}(B-p+2) \cdots (t-T)^{-1}
(B-1)\cdot (t-T)^{-1}B\mod \p.$$ \Endproof
\end{lem}

\begin{rem}\label{remput}{\rm The above proposition is interesting because there is no 
closed formula known for the computation of the $\p$-curvature 
if the rank is $>2$ (compare
to \cite{vdutComp} 
and \cite{vdPutpcurv}). The above lemma yields such 
a closed formula for Okubo systems. On the other hand, every irreducible 
Fuchsian system is a subfactor of an Okubo system (this is known 
by the work of Okubo and
follows also
from Prop. \ref{proptriv} and the fact that, by Equation \eqref{eq211} below, 
 the middle convolution 
$\mc_{-1}$ induces  the identity tranformation of Fuchsian systems).} 
\end{rem}

The following technical proposition will be used below:

\begin{prop}\label{nilpo} Let $K$ be 
a number field, let $F:Y'=CY,\,C\in K(t)^{n\times n},$ be a
 Fuchsian system 
 and let 
 $\c_0(F):Y'=DY$ be the naive convolution of $F$ with the parameter $0.$  
Let $p$ be a prime number 
and let 
$\p$  denote a prime of $K$ which lies
over $p$ such that  the entries of $C$ can be 
reduced modulo $\p.$ If 
$C_p^k \equiv 0 \mod \p,$ then  
 $D_{kp+1}\equiv 0 \mod \p.$ 
\end{prop}

\proof The Fuchsian system $F$ is the Fuchsian system associated to 
some tuple $(A_1,\ldots,A_r)\in (K^{n\times n})^r.$ By Prop.~\ref{proptriv},
the naive convolution of $F$ with the parameter $-1$ can be written 
in Okubo normal form
$$\c_{-1}(F):Y'=(t-T)^{-1}\tilde{B}Y=\tilde{D}Y\,.$$
Here, $T$ is the diagonal matrix 
$T={\rm diag}(t_1,\ldots,t_1,\ldots,t_r,\ldots,t_r)$
(every $t_k$ occurs $n$ times) and 
  $$\tilde{B}=(\tilde{B}_{i,j}),\; \tilde{B}_{i,j}=A_j -\delta_{i,j}E_n$$ 
($E_n\in \GL_n(K)$ denoting the identity matrix). 
  Using the  base change  which is induced by 
 the matrix  $H(t-T),$ where
 \[H:=\left( \begin{array}{cccc}
       E_n & -E_n & 0 & \ldots \\
       0  &      \ddots&\ddots    \\
       \vdots&  & &- E_n \\
       0  & \ldots&& E_n
    \end{array} \right), \]
one can verify that the naive convolution 
 $\c_{-1}(F)$ is equivalent to the following 
system:
 \begin{equation}\label{eq211}  
  Y'=GY= \left( \begin{array}{cccc}
       0 &  \ldots & \ldots & 0 \\
            \vdots &&&\vdots\\
                         0 & \ldots & \ldots &0 \\
      {{A_1}\over{t-t_1}} &({{A_1}\over{t-t_1}}+{{A_2}\over{t-t_2}}) 
&\ldots & ({A_{1}\over{t-t_1}}+\cdots + 
{A_{r}\over{t-t_r}})
    \end{array} \right)Y \,.\end{equation} By a straightforward computation, one sees that 
\begin{equation}\label{eq21}
G = H(t-T)\cdot  D \cdot (H(t-T))^{-1}\,, \end{equation}
where $D$ is the matrix appearing in 
the naive convolution $\c_0(F):Y'=DY$ with the parameter $0.$
Using the Leibnitz rule, one can further  verify that
\begin{equation}\label{eq22} G_p= (H(t-T))\cdot  \tilde{D}_p\cdot  (H(t-T))^{-1}\,\end{equation} (note that $\tilde{D}_p$ belongs to the 
naive convolution $c_{-1}(F)$ with parameter $-1$).
It follows from Equation \eqref{eq211} that 
 \[ G_p= \left( \begin{array}{cccc}
       0 &  \ldots & \ldots &0 \\
            \vdots &&&\vdots\\
                         0 & \ldots&\ldots  & 0 \\
      \ast & \ast &\ast  & C_p
    \end{array} \right) \,. \]
Using this and the assumption  $C_p^k\equiv  0\mod \p$ one sees that 
\[ G_p^k G = \left( \begin{array}{cccc}
       0 &  \ldots & \ldots &0 \\
            \vdots &&&\vdots\\
                         0 & \ldots&\ldots  & 0 \\
      \ast & \ast &\ast  & 0
    \end{array} \right) \left( \begin{array}{cccc}
       0 &  \ldots & \ldots &0 \\
            \vdots &&&\vdots\\
                         0 & \ldots&\ldots  & 0 \\
      \ast & \ast &\ast  & \ast
    \end{array} \right)   \equiv 0\mod \p    \,. \]
Thus, by Equations \eqref{eq21} and \eqref{eq22}, 
\[H(t-T)\cdot 
{\tilde{D}_p}^k \cdot  D\cdot 
(H(t-T))^{-1}\equiv 0 \mod \p\]
and thus 
$$ {\tilde{D}_p}^k \cdot  {D}\equiv 0 \mod \p\,.$$
 The claim 
follows now from this and the equality 
 \[ {\tilde{D}_p}^k\cdot  {D}\equiv D_{pk+1} \mod \p \,,\]
which is immediate  from the Equation \ref{eq2111} appearing in
Lemma \ref{recu}.
\Endproof

\begin{thm}\label{cor21} Let $K$ be 
a number field, let $$F:Y'=CY,\,C\in K(t)^{n\times n},$$ be a
 Fuchsian system and let $\mu \in \QQ.$ 
Let $p$ be a prime number 
and let 
$\p$  denote a prime of $K$ which lies
over $p$ such that the number $\mu$ and 
 the entries of $C$ can be 
reduced modulo $\p.$ If the $\p$-curvature of $F$ is 
nilpotent of rank $k$ (i.e.,
$\bar{C}_p(F)^k \equiv 0$),
then  the $\p$-curvature of the middle convolution 
$\mc_\mu(F)$ is nilpotent of 
rank $r\in \{k-1,k,k+1\}.$
\end{thm}

\proof We first prove that the $\p$-curvature of the naive convolution
$\c_\mu(F)$ is nilpotent of rank at most $k+1:$
Let 
 $\c_0(F):Y'=DY$ be the naive convolution of $F$ with the parameter $0$
and let  $\c_\mu(F):Y'=D^\mu Y$ 
be the naive convolution of $F$ with the parameter $\mu.$
 By assumption, one has 
$C_p^k \equiv 0 \mod \p\,.$ Thus, by the preceding 
proposition, 
 \begin{equation}\label{eq23}D_{kp+1}\equiv 0 \mod \p\,.\end{equation}
One can easily deduce  from  the Equation \eqref{eq2111} of 
 Lemma \ref{recu}, that 
the matrix $D_{kp+1}\mod \p$ 
 appears as a factor in the product formula \eqref{eq2111} 
for $D^\mu_{p(k+1)}$ taken modulo $\p.$ 
Thus, by Equation \eqref{eq23} and again by Lemma \ref{recu}, 
 \[0\equiv D^\mu_{p(k+1)}  \equiv 
 (D^\mu_p)^{k+1}=\bar{C}_p(\mc_\mu(F))^{k+1}\mod \p \,.\]

Since the middle convolution $\mc_\mu(F)$
is a factor of the naive convolution
$\c_\mu(F),$ it follows that the $\p$-curvature of the middle convolution 
$\mc_\mu(F)$ is nilpotent of 
rank at most $k+1.$ Now the claim follows from the fact that 
$\mc_{-\mu}\circ \mc_{\mu}=\id.$
\Endproof

\section{Two new examples of globally nilpotent Fuchsian systems}
\label{Section4}

\subsection{A globally nilpotent Fuchsian system of rank two.}\label{sec41}

Dwork has conjectured that any globally nilpotent second 
order differential equation has either algebraic solutions or has a 
correspondence to a Gau{\ss} hypergeometric differential equation (see
\cite{Krammer}).
This conjecture was disproved by Krammer (loc.~cit.). The 
counterexample which  Krammer
gives is the uniformizing 
differential equation $K$ of an arithmetic 
 Fuchsian lattice. It 
comes from the periods of a family of abelian surfaces  over 
a Shimura curve $S\simeq \PP^1\setminus \{0,1,81,\infty\}.$ The 
equation $K$  is an 
irreducible 
ordinary  second order 
differential equation which can easily be transformed into the following 
 Fuchsian 
system of rank two: 
$$ K:Y'
=\left(\, \frac{1}{t}\cdot \left(\begin{array}{cc} 0&0\\
-\frac{1}{2}& -\frac{1}{2} \end{array}\right)+ 
\frac{1}{t-1}\cdot\left(\begin{array}{cc} 0&0\\
\frac{4}{9}& -\frac{1}{2} \end{array}\right)+
\frac{1}{t-81}\left(\begin{array}{cc} 0&1\\
0& \frac{1}{2} \end{array}\right)\,
\right)\cdot Y .$$
The local monodromy of 
$K$ at the finite singularities is given by three reflections,
the local monodromy of 
$K$  at $\infty$ is given by  an element of order six. \\

\begin{thm} \label{thm41}
The middle convolution $H:=\mc_{\frac{1}{6}}(K)$
is equivalent to the following Fuchsian system:
$$Y'=
 \left(\,\frac{1}{t} \left(\begin{array}{cc} -\frac{19}{30}&\frac{19}{10}\\
-\frac{1}{10}& \frac{3}{10} \end{array}\right)
 + 
\frac{1}{t-1}\left(\begin{array}{cc} -\frac{1}{3}&-\frac{7}{18}\\
0&0 \end{array}\right)+
\frac{1}{t-81}\left(\begin{array}{cc} 0&0\\
-\frac{1}{2} & \frac{2}{3} \end{array}\right)\,
\right)Y  .$$
The Fuchsian system $H$ is 
 globally nilpotent. Moreover, the 
system $H$ has neither a correspondence to a hypergeometric 
differential system, nor a correspondence to a uniformizing 
differential equation of a Fuchsian lattice.
\end{thm}

\proof 
The first claim follows from Cor. \ref{cor21}. To
prove the last statement, we use the Riemann-Hilbert correspondence:
The monodromy tuple of $K$ can be shown to be 
$$ A:=\left( i\cdot 
 \left(\begin{array}{cc}
  0&-(\wV+1)/2 \\
    (\wV-1)/2 &0
  \end{array}\right),\,\,
i\cdot \left(\begin{array}{cc}
   -\wIII &1+\wV\\
    1-\wV& \wIII 
  \end{array}\right) \right., $$
$$
 \left. i\cdot \left(\begin{array}{cc}
   (\wIII+\sqrt{15})/2& -2 \\
   2&(\wIII-\sqrt{15})/2
  \end{array}\right)
\right)\,.$$ In the following, the element $\zeta_n$ denotes 
the root of unity $e^{2\pi i /n}.$
Using  Thm. \ref{monodromy}, it is easy to see that the monodromy of
$\mc_{\frac{1}{6}}$ is given by 
$$ \MC_{\zeta_6}(A)=(B_1,B_2,B_3)=$$

$$\left(
 \zeta_6\cdot  \left( \begin{array}{cc}
   \frac{1}{5}(\wV + 3\zeta_{60}^{10} + 6\zeta_{60}^8   - 
    6\zeta_{60}^2 - 3)&
      \frac{1}{5} (-6 \wV - 10\zeta_{60}^{10} - 10\zeta_{60}^8  + 10\zeta_{60}^2 + 10) \\
 \frac{1}{5}(\zeta_{60}^{14} + 2\zeta_{60}^{10} - \zeta_{60}^8 - \zeta_{60}^6 - \zeta_{60}^4 + \zeta_{60}^2 
    + 1) & \frac{1}{5}(-\wV - 3\zeta_{60}^{10} - 6\zeta_{60}^8)
   \end{array}\right)\right.,$$
$$ 
\left.\zeta_6\cdot \left(\begin{array}{cc}
   -\zeta_{60}^{10} &  2\zeta_{60}^{10} \\
             0 & \zeta_{60}^{10} - 1
\end{array}\right),
\zeta_6\cdot\left(\begin{array}{cc}
 \zeta_{60}^{10} - 1  &             0 \\
-\zeta_{60}^{10} + 1  &   -\zeta_{60}^{10}
\end{array}\right) \right).$$ 
By \cite{DR}, Cor 5.9, one can show that 
the elements $B_i,\,i=1,2,3,$ generate a subgroup which is conjugate 
to a subgroup of 
${\rm SU}_{1,1}(\RR).$ The latter group is well known 
to be  conjugate to the group
$\SL_2(\RR)$ inside the group $\GL_2(\CC).$

Consider the element 
$$\tilde{B}:=B_1B_2B_1=$$
$$\left(
\begin{array}{cc}
 2\zeta_{60}^{10} - 1 & -4\zeta_{60}^{10} - 2\zeta_{60}^8 + 2\zeta_{60}^2\\
-\zeta_{60}^{14} + \zeta_{60}^8 + \zeta_{60}^6 + \zeta_{60}^4 - \zeta_{60}^2 - 1
 & 2\zeta_{60}^{14} - 
    2\zeta_{60}^{10} - 2\zeta_{60}^6 - 2\zeta_{60}^4 + 3
\end{array}
\right)\,.$$  It is an elliptic element 
since the absolute value of its trace is $|-\wV+1|<2.$
 Moreover, the order of $\tilde{B}$ is 
infinite. (The eigenvalues of $\tilde{B}$ are roots of
 $f:=X^4-2X^3-2X^2-2X+1$ with $\Gal(f)=D_8,$ thus the eigenvalues 
cannot be  roots of unity.)
 Thus, the monodromy tuple  $\MC_{\zeta_6}(A)$ of 
$\mc_{\frac{1}{6}}(F)$ generates a non-discrete subgroup of $\PGL_2(\RR)$ and 
thus does not have any correspondence to a hypergeometric 
differential system (in which case it would be a 
discrete triangle group, see \cite{Doran}, 4.1.2),
 nor a correspondence to a uniformizing 
differential equation of a Fuchsian lattice.
\Endproof

\subsection{A globally nilpotent system with 
differential Galois group $G_2$.}\label{sec42}

The following construction of local systems
with finite monodromy will be used in Section 
\ref{sec42} below: Let $\rho_f:\pi_1(X)\to G$ be 
the surjective homomorphism associated to a finite 
Galois cover with Galois group $G$ 
and let $\alpha:G\to \GL_n(\CC)$ be a representation.
Then the composition 
$$ \alpha\circ \rho_f: \pi_1(X)\to \GL_n(\CC)$$ corresponds to a 
local system
 with finite monodromy on $X$ which is denoted by $\V_{f,\alpha}.$ 

Consider the following local systems 
 with finite monodromy $$\V_i\; =\; \V_{f_i,\alpha},\quad i=0,1,2,$$
on $S =\AA^1\setminus \{0,1\}$,  where the 
$f_i$ and $\alpha$ are as follows: The representation  
$\alpha$ is the embedding of $\ZZ/2\ZZ$ into ${\QQ_\ell}^\times$ 
and  the $f_i:Y_i \to S,\, i=0,1,2,$ 
are the    double covers defined by the equations
$$ z^2=t(t-1),\quad {\rm resp.}\quad z^2=t-1,\quad {\rm resp.}\quad z^2=t\,.$$
The middle convolution operation we use is $\MC_\chi,$ where 
$\chi$ is associated to the double cover of $\GG_m$ 
given by $z^2=t$ and by $\alpha:\ZZ/2\ZZ\to \QQ_\ell$ as above. 
Let $\GGG$ denote the local system
 on $S,$ defined by 
\begin{equation}\label{eq1} 
 \V_2\otimes(\MC_\chi(\V_1\otimes \MC_\chi(\V_2\otimes\MC_\chi(\V_1\otimes \MC_\chi(\V_2\otimes \MC_\chi(\V_1\otimes\MC_\chi(\V_0))))))))\,.\end{equation}

In \cite{DR06}, the following is shown:

\begin{prop}\label{prop41} The monodromy tuple of $\GGG$ is of the form 
$(\rho_{\GGG}(\delta_1),\rho_{\GGG}(\delta_2))\in \GL_7(\ZZ),$
 where $(\rho_{\GGG}(\delta_1),\rho_{\GGG}(\delta_2))$ is as follows:
\begin{small}
$$\left(\left(\begin{array}{ccccccc}
 1&  0&  0&  2&  2&  0&  0\\
 0 & 1 & 0& -2 & 0&  2 & 0\\
 0 & 0 & 1&  2&  2&  2 & 2\\
 0 & 0&  0& -1&  0&  0&  0\\
 0 & 0&  0 & 0 &-1 & 0 & 0\\
0  &0 & 0 & 0&  0& -1 & 0\\
 0&  0 & 0&  0 & 0&  0& -1
\end{array}
\right),\quad
\left(\begin{array}{ccccccc}
1 &0& 0& 0& 0& 0& 0\\
0& 1& 0& 0& 0& 0& 0\\
0& 0& 1& 0& 0& 0& 0\\
2& 0& 0& 1 &0& 0& 0\\
0& 2& 0& 0& 1& 0& 0\\
0& 0& 2& 0& 0& 1& 0\\
0& 0& 0 &2 &4 &4& 1
\end{array}
\right)\right)\,.$$\end{small}
For any prime $\ell>3,$ the $\ell$-adic 
closure of the monodromy group coincides with the group 
$G_2(\ZZ_\ell).$ Especially, the Zariski closure of the monodromy 
group in $GL_7(\CC)$
coincides with the group $G_2(\CC).$
\end{prop}

Now, consider the following Fuchsian rank one systems $F_i,\,i=0,1,2:$ 
$$ F_0:Y'=\left(-\frac{1}{2t}-\frac{1}{2(t-1)}\right) Y,\quad 
F_1:Y'=\frac{1}{2(t-1)} Y =\left(\frac{0}{2t} +\frac{1}{2(t-1)}\right) Y$$
and  
$$ F_2:Y'=\frac{1}{2t} Y= \left(\frac{1}{2t} +\frac{0}{2(t-1)}\right) Y \,.$$ 
The local systems of holomorphic solutions of 
$F_i,\, i=0,1,2,$ on $S=\AA^1\setminus \{0,1\}$ 
can easily be seen to be isomorphic to 
the above local systems $\V_i,\, i=0,1,2.$\\

The tensor 
product of two Fuchsian systems associated to 
$(A_1,\ldots, A_r)$ and $(B_1,\ldots,B_r)$ 
 is given by the Fuchsian system associated to the tuple
$$\left(A_1\otimes 1+1\otimes B_1,\ldots,A_r\otimes 1+1\otimes  B_r\right)\,.$$
The dual of a system $F:Y'=CY$ of linear differential equations 
is given by $F^*:Y'=-C^{\rm tr}Y,$ where $C^{\rm tr}$ denotes the transpose of
the matrix $C.$\\

 Let $G$ be the 
Fuchsian system which is given by the following sequence of 
middle convolutions and tensor products:
\begin{equation}\label{eq2} 
 F_2^*\otimes(\mc_{-\frac{1}{2}}(F_1^*\otimes \mc_{\frac{1}{2}}(F_2^*\otimes\mc_{-\frac{1}{2}}(F_1\otimes \mc_{\frac{1}{2}}(F_2\otimes \mc_{-\frac{1}{2}}
(F_1^*\otimes\mc_{\frac{1}{2}}(F_0))))))))\,.\end{equation}
 
Recall from \cite{SingervanderPut}, that 
any  system of linear 
differential equations $$F:Y'=AY, \,\, A \in \CC(t)^{n\times n},$$ 
 has attached
an algebraic  group $G_F$ to it, called the {\em differential Galois group}
of $F.$ Moreover, if the singularities of $F$ are regular, then 
$G_F$ can be viewed as the Zariski closure of the monodromy of $F$ 
inside
the group $\GL_n(\CC).$ With this, we obtain the following result:

\begin{thm}\label{thm421} The Fuchsian system $G$ is given by
\begin{tiny}
\[ Y'=\left( \frac{1}{2t}\cdot \left(\begin{array}{ccccccc}
   -2  &  0  &  0  &  0  &  0  &  0  &  0\\
    0 &  -2  &  0  &  0  &  0  &  0  &  0\\
    0  &  0 &  -2  &  0  &  0  &  0  &  0\\
    0  &  0 &-1 &-1  &  0  &  0  &  0\\
-1  &  0  &  0  &  0 &-1  &  0  &  0\\
    0 &-1  &  0  &  0  &  0 &-1  &  0\\
    0  &  0& -1  &  0  &  0  &  0 &-1\\
\end{array}\right)+
 \frac{1}{2(t-1)}\cdot \left(\begin{array}{ccccccc}
    0  &  0  &  0  &  0 &  -2 & 1 & 1\\
    0  &  0  &  0  &  0  &  0  &  0 &-1\\
    0  &  0  &  0  &  0 &-1  &  0  &  0\\
    0  &  0  &  0 &  -2  &  0  &  0  &  0\\
    0  &  0  &  0  &  0 &  -2  &  0  &  0\\
    0  &  0  &  0  &  0  &  0 &  -2  &  0\\
    0  &  0  &  0 & -1 & 1 &-1 &  -2\\
\end{array}\right)\right)Y. \]\end{tiny}
 The differential Galois group of 
$G$ is isomorphic to $G_2(\CC).$ Moreover, the Fuchsian system $G$ 
 is  globally nilpotent.
\end{thm}
\proof The form of the matrices follows from the explicit construction 
of $\mc_\mu.$ By the Riemann-Hilbert correspondence for $\mc_\mu$ (Thm. \ref{monodromy}),
the monodromy tuple of $G$ is given by the matrices 
$(\rho_{\GGG}(\delta_1),\rho_{\GGG}(\delta_2))$ above. Thus the first claim 
follows from Prop. \ref{prop41} and the fact that 
the differential Galois group of a Fuchsian system is the Zariski closure 
of the monodromy (since all 
singularities are regular), see \cite{SingervanderPut}. 
The last claim follows from Thm. \ref{cor21}.
\Endproof
 
\begin{rem}{\rm The above Fuchsian system can be viewed to be 
the de Rham realization of a certain mixed motive $\M$ of dimension 
$7,$ which appears in \cite{DR06}. The motive $\M$ stands in close connection 
to a question of Serre on the existence 
of motives with exceptional Galois group $G_2,$ see \cite{Ser2}, \cite{DR06}.}
\end{rem}

\bibliographystyle{plain} 

\bibliography{p}

Michael Dettweiler

IWR, Universit\"at Heidelberg,

INF 368

69121 Heidelberg, Deutschland

e-mail: michael.dettweiler@iwr.uni-heidelberg.de\\

Stefan Reiter

Technische Universit\"at Darmstadt

Fachbereich Mathematik AG 2

Schlo\ss gartenstr. 7

64289 Darmstadt, Deutschland

email: reiter@mathematik.tu-darmstadt.de

\end{document}